\documentclass[12pt]{article}
\usepackage{amsfonts}
\usepackage{amsmath}
\usepackage{amssymb}

\newtheorem{theorem}{Theorem}[section]

\newtheorem{example}[theorem]{Example}

\newtheorem{proposition}[theorem]{Proposition}

\input{tcilatex}

\begin{document}

\title{The $\mathfrak{n}$-homology of representations. }
\author{Tim Bratten \\
Facultad de Ciencias Exactas, UNICEN. Tandil, Argentina. }
\date{}
\maketitle

\begin{abstract}
The $\mathfrak{n}$-homology groups of a $\mathfrak{g}$-module provide a
natural and fruitful extension of the concept of highest weight to the
representation theory of a noncompact reductive Lie group. In this article
we give an introduction to the $\mathfrak{n}$-homology groups and survey
some developments, with a particular emphasis on results pertaining to the
problem of caculating $\mathfrak{n}$-homology groups.
\end{abstract}

\section{Introduction}

The concept of a highest weight and its use to classify irreducible
representations of compact Lie groups can traced back nearly a century, to
seminal work by E. Cartan and H. Weyl. For a compact, connected Lie group,
the highest weight theory gives a tight parametrization of irreducible
representations in terms of specific invariants associated to the group. If
one tries to extend this concept to the representation theory of a
noncompact, real reductive group one immediately encounters two problems. On
the one hand, in the noncompact case, it turns out there are several
conjugacy classes of complex Borel subalgebras, and what might be called a
highest weight depends on the choice of a conjugacy class. On the other
hand, it is quite common that what should be called a highest weight turns
out to be zero for every choice of Borel subalgebra. This means, in the
traditional sense, the highest weight does not exist for a great majority of
irreducible representations.

Although there is no way to avoid the first problem, representation
theorists have confronted the second problem by considering the highest
weight to be a functorial construction and studying the related derived
functors. This has proved to be especially fruitful, producing a strong and
useful family of invariants associated to a representation. In this article
we give a brief introduction to the $\mathfrak{n}$-homology (and $\mathfrak{n%
}$-cohomology) groups, followed by a survey of some results, focusing on
developments related to the problem of calculating the $\mathfrak{n}$%
-homology of representations.

The author would like to thank the organizers of the 2007 meeting of the UMA
for giving him an opportunity to present some results in the form of a
conference and for asking him to submit this article. He would also like to
acknowledge the help and encouragement he has received from Jorge Vargas.
This article is dedicated to the memory of Misha Cotlar, with a special
recognition to Dr. Cotlar's role as advisor and mentor to the late Jos\'{e}
Pererya.

\section{$\mathfrak{n}$-homology and $\mathfrak{n}$-cohomology}

In this section we introduce the $\mathfrak{n}$-homology and $\mathfrak{n}$%
-cohomology of $\mathfrak{g}$-modules (for more details see \cite{KV}).

Let $\mathfrak{g}$ be a complex reductive Lie algebra. By definition, \emph{%
a Borel subalgebra }of $\mathfrak{g}$ is a maximal solvable subalgebra and 
\emph{a parabolic subalgebra }of $\mathfrak{g}$ is a subalgebra that
contains a Borel subalgebra. If $\mathfrak{p}\subseteq \mathfrak{g}$ is a
parabolic subalgebra then \emph{the nilradical} $\mathfrak{n}$ of $\mathfrak{%
p}$ is the largest solvable ideal in $\left[ \mathfrak{p},\mathfrak{p}\right]
$. \emph{A\ Levi factor }is a complementary subalgebra to $\mathfrak{n}$ in $%
\mathfrak{p}$. One knows that Levi factors exist and that they are exactly
the subalgebras which are maximal with respect to being reductive in $%
\mathfrak{p}$. When $\mathfrak{l}$ is a Levi factor than 
\begin{equation*}
\mathfrak{p=l}\oplus \mathfrak{n}
\end{equation*}%
is called \emph{a Levi decomposition}.

Fix a parabolic subalgebra $\mathfrak{p}$ with nilradical $\mathfrak{n}$ and
Levi factor $\mathfrak{l}$. Let $U(\mathfrak{n})$ denote the enveloping
algebra of $\mathfrak{n}$ and let $\mathbb{C}$ be the 1-dimensional trivial
module. If $M$ is a $\mathfrak{g}$-module then \emph{the zero} $\mathfrak{n}%
\emph{-}$\emph{homology of} $M$ is the $\mathfrak{l}$-module 
\begin{equation*}
H_{\text{0}}(\mathfrak{n},M)=\mathbb{C}\otimes _{U(\mathfrak{n})}M.
\end{equation*}%
This $\mathfrak{l}$-module is sometimes referred to as \emph{the space of
coinvariants}, although it clearly depends on the choice of parabolic
subalgebra. The definition of the zero homology determines a right exact
functor from the category of $\mathfrak{g}$-modules to the category of $%
\mathfrak{l}$-modules. \emph{The} $\mathfrak{n}$-\emph{homology groups of} $%
M $ are the $\mathfrak{l}$-modules obtained as the corresponding derived
functors. There is a standard complex for calculating these homology groups,
defined as follows. \emph{The right standard resolution of }$\mathbb{C}$ is
the complex of free right $U(\mathfrak{n})$-modules given by 
\begin{equation*}
\cdots \rightarrow \Lambda ^{p+1}\mathfrak{n}\otimes U(\mathfrak{n}%
)\rightarrow \Lambda ^{p}\mathfrak{n}\otimes U(\mathfrak{n})\rightarrow
\cdots \rightarrow \mathfrak{n}\otimes U(\mathfrak{n})\rightarrow U(%
\mathfrak{n})\rightarrow 0\text{.}
\end{equation*}%
Applying the functor 
\begin{equation*}
-\otimes _{U(\mathfrak{n})}M\text{ }
\end{equation*}%
to the standard resolution we obtain a complex 
\begin{equation*}
\cdots \rightarrow \Lambda ^{p+1}\mathfrak{n}\otimes M\rightarrow \Lambda
^{p}\mathfrak{n}\otimes M\rightarrow \cdots \rightarrow \mathfrak{n}\otimes
M\rightarrow M\rightarrow 0
\end{equation*}%
of left $\mathfrak{l}$-modules called \emph{the standard} $\mathfrak{n}$%
\emph{-homology complex}. Here $\mathfrak{l}$ acts via the tensor product of
the adjoint action on $\Lambda ^{p}\mathfrak{n}$ with the given action on $M$%
. Since $U(\mathfrak{g})$ is a free $U(\mathfrak{n})$-module, a routine
homological argument identifies the pth homology of the standard complex
with the pth $\mathfrak{n}$-homology group 
\begin{equation*}
H_{\text{p}}(\mathfrak{n},M).
\end{equation*}%
One can prove that the induced $\mathfrak{l}$-action on the homology groups
of the standard complex is the correct one.

\emph{The zero} $\mathfrak{n}$\emph{-cohomology }of a $\mathfrak{g}$-module $%
M$ is the $\mathfrak{l}$-module 
\begin{equation*}
H^{0}(\mathfrak{n},M)=\text{Hom}_{U(\mathfrak{n})}(\mathbb{C},M).
\end{equation*}%
This $\mathfrak{l}$-module is sometimes referred to as \emph{the space of
invariants}, and also clearly depends on the choice of parabolic subalgebra
The definition of the zero cohomology determines a left exact functor from
the category of $\mathfrak{g}$-modules to the category of $\mathfrak{l}$%
-modules. By definition, \emph{the} $\mathfrak{n}$\emph{-cohomology groups of%
} $M$ are the $\mathfrak{l}$-modules obtained as the corresponding derived
functors. These $\mathfrak{l}$-modules can be calculated by applying the
functor 
\begin{equation*}
\text{Hom}_{U(\mathfrak{n})}(-,M)
\end{equation*}%
to the standard resolution of $\mathbb{C}$, this time by free left $U(%
\mathfrak{n})$-modules. In a natural way, one obtains a complex of $%
\mathfrak{l}$-modules and the pth cohomology of this complex realizes the
pth $\mathfrak{n}$-cohomology group

\begin{equation*}
H^{\text{p}}(\mathfrak{n},M).
\end{equation*}

It turns out that the structure of the $\mathfrak{n}$-cohomology is
determined by the structure of the $\mathfrak{n}$-homology, in a simple way
. Thus, it is often a matter of convenience whether one works with homology
groups or cohomology groups. In this article, we will focus on results
framed in terms of homology. The following proposition, whose proof is
established by an analysis of standard complexes, can be used to translate
results about $\mathfrak{n}$-homology into results about $\mathfrak{n}$%
-cohomology \cite[Section 2]{HS}.

\begin{proposition}
Suppose $M$ is a $\mathfrak{g}$-module. Let $\mathfrak{p}\subseteq \mathfrak{%
g}$ be a parabolic subalgebra with nilradical $\mathfrak{n}$ and Levi factor 
$\mathfrak{l}$. Let d denote the dimension of $\mathfrak{n}$. Then there are
natural isomorphisms 
\begin{equation*}
H_{\text{p}}(\mathfrak{n},M)\cong H^{\text{d-p}}(\mathfrak{n},M)\otimes
\Lambda ^{\text{d}}\mathfrak{n}\newline
.
\end{equation*}
\end{proposition}

\section{Representations of linear reductive Lie groups}

In this section we review some classical results about the representation
theory of reductive Lie groups (for details see \cite{Wa}) and introduce the
canonical globalizations.

For simplicity we work with a class of reductive Lie groups we call \emph{%
linear}, although there is no problem working in the more general context of
a reductive group of Harish-Chandra class. In particular, we assume the
following setup. $G$ will denote a connected, complex reductive group. This
means $G$ is a connected, complex Lie group with the property that the
maximal compact subgroups are real forms of $G$. The group $G_{0}$ will
denote a real form of $G$ and is assumed to have finitely many connected
components. We call $G_{0}$ \emph{a linear reductive Lie group.} The Lie
algebras of $G$ and $G_{0}$ will be denoted $\mathfrak{g}$ and $\mathfrak{g}%
_{0}$, respectively. For the remainder of this article we fix a maximal
compact subgroup $K_{0}$ of $G_{0}$ and let $K\subseteq G$ be the
complexification of $K_{0}$. In general, we write $K,L$ etc. to indicate
complex subgroups of $G$ and denote the corresponding Lie algebras by $%
\mathfrak{k},\mathfrak{l}$ etc. Subgroups of $G_{0}$ will be denoted by $%
K_{0},L_{0}$ etc. with the corresponding real Lie algebras written as $%
\mathfrak{k}_{0},\mathfrak{l}_{0}$ etc.

A \emph{representation} of $G_{0}$ will mean a continuous linear action of $%
G_{0}$ in a complete, locally convex topological vector space. When we speak
of irreducible or finite length representations, the corresponding
definitions should be framed in terms of invariant closed subspaces. A
vector $v$ in a representation $V$ is called \emph{smooth} when 
\begin{equation*}
\lim_{t\rightarrow 0}\frac{\exp (t\xi )v-v}{t}\text{ \ exists for each }\xi
\in \mathfrak{g}_{0}.
\end{equation*}%
In order to define $\mathfrak{n}$-homology groups, we will be primarily
interested in \emph{smooth representations}. These are representations where
every vector is smooth. In a natural way a smooth representation carries a
compatible $\mathfrak{g}$-action. For a compact Lie group one can show that
a finite length representation is finite-dimensional and therefore smooth.

We recall some basic results about the infinite-dimensional representations
of reductive groups. In the 1950s, Harish-Chandra proved that an irreducible
unitary representation $V$ has the property that each irreducible $K_{0}$%
-submodule has finite multiplicity in $V$. This led him to define and study 
\emph{admissible representations}. This means that each irreducible $K_{0}$%
-submodule of the representation has finite multiplicity. Harish-Chandra
then considered \emph{the subspace of }$K_{0}$-\emph{finite vectors. }By
definition, a vector $v$ in a representation is called $K_{0}$-\emph{finite}
if the span of the $K_{0}$-orbit of $v$ is finite-dimensional. Although the
subspace of $K_{0}$-finite vectors is not $G_{0}$-invariant, Harish-Chandra
proved that $K_{0}$-finite vectors are smooth, and thus form a $(\mathfrak{g}%
,K_{0})$-module called \emph{the underlying Harish-Chandra module.}

On the other hand, it is possible to define abstractly the concept of a
Harish-Chandra module. This is a $\mathfrak{g}$-module equipped with a
compatible, locally finite $K_{0}$-action. Harish-Chandra proved that an
irreducible Harish-Chandra module appears as the underlying $(\mathfrak{g}%
,K_{0})$-module of $K_{0}$-finite vectors in an irreducible admissible
Banach space representation for $G_{0}$ and W. Casselman proved that the
same holds for any finite-length Harish-Chandra module. By now we know more.
In particular, given a Harish-Chandra module $M$ we define a globalization $%
M_{\text{glob}}$ of $M$ to be an admissible representation for $G_{0}$ whose
underlying Harish-Chandra is $M$. We assume our Harish-Chandra modules have
finite-length. Then we can assert that several canonical and functorial
globalizations exist on the category of Harish-Chandra modules. These are:
the smooth globalization of Casselman and Wallach \cite{C}, its dual
(called: the distribution globalization), Schmid's minimal globalization 
\cite{S} and its dual (the maximal globalization). All four globalizations
are smooth. We will let $M_{\text{min}}$, $M_{\text{max}}$, $M_{\infty }$
and $M_{\text{dis}}$ denote respectively, the minimal, the maximal, the
smooth and the distribution globalizations of a Harish-Chandra module $M$.
If $M_{\text{glob}}$ denotes a Banach globalization of $M$, then there is a
natural chain of inclusions 
\begin{equation*}
M\subseteq M_{\text{min}}\subseteq M_{\infty }\subseteq M_{\text{glob}%
}\subseteq M_{\text{dis}}\subseteq M_{\text{max}}.
\end{equation*}%
\noindent In this chain the minimal globalization\emph{\ }is known to
coincide with the analytic vectors in $M_{\text{glob}}$ while $M_{\infty }$
coincides with the smooth vectors in $M_{\text{glob}}$. In particular, one
knows that a finite-length admissible Banach space representation for $G_{0}$
is smooth if and only if it is finite dimensional. Later in this article we
will review various results, often called \emph{comparison theorems},
relating the $\mathfrak{n}$-homologies of a Harish-Chandra module to the $%
\mathfrak{n}$-homologies of a canonical globalization.

\section{Some structural details}

In this section we recall some structure theory and an important technical
result about the decomposition of $\mathfrak{n}$-homology groups for certain 
$\mathfrak{g}$-modules, giving special emphasis on the case of a Borel
subalgebra. Recall that a Cartan subalgebra $\mathfrak{h}\subseteq \mathfrak{%
g}$ is a maximal abelian subalgebra whose elements are semisimple under the
adjoint representation of $\mathfrak{h}$ in $\mathfrak{g}$ . A nonzero
eigenvalue $\alpha \in \mathfrak{h}^{\ast }$ for the adjoint representation
is called a root. $\Sigma $ will denote the set of roots. Thus 
\begin{equation*}
\mathfrak{g}=\mathfrak{h\oplus }\dbigoplus\limits_{\alpha \in \Sigma }%
\mathfrak{g}^{\alpha }\text{ \ }
\end{equation*}%
where $\mathfrak{g}^{\alpha }$ is the eigenspace corresponding to root $%
\alpha $. One knows that $\alpha \in \Sigma $ if and only if $-\alpha \in
\Sigma $. If $\mathfrak{b}$ is a Borel subalgebra of $\mathfrak{g}$
containing $\mathfrak{h}$ then the roots of $\mathfrak{h}$ in $\mathfrak{b}$
define a subset $\Sigma ^{+}\subseteq \Sigma $ called \emph{the
corresponding set of positive roots}. When the sum of two positive roots is
a root, then that sum is positive. One also knows that $\Sigma $ is a
disjoint union: 
\begin{equation*}
\Sigma =\Sigma ^{+}\cup -\Sigma ^{+}\text{.}
\end{equation*}%
One can show there is a unique $H_{\alpha }\in \left[ \mathfrak{g}^{\alpha },%
\mathfrak{g}^{-\alpha }\right] $ such that $\alpha \left( H_{\alpha }\right)
=2.$We use this element to define the value of \emph{the dual root}. In
particular, the dual root is given by 
\begin{equation*}
\overset{\vee }{\alpha }(\mu )=\mu (H_{\alpha })\text{ \ for }\mu \in 
\mathfrak{h}^{\ast }\text{.}
\end{equation*}%
The linear reflection $s_{\alpha }:\mathfrak{h}^{\ast }\rightarrow \mathfrak{%
h}^{\ast }$ corresponding to $\alpha \in \Sigma $ is defined as 
\begin{equation*}
s_{\alpha }(\mu )=\mu -\overset{\vee }{\alpha }(\mu )\alpha .
\end{equation*}%
These reflections generate a finite subgroup of the general linear group of $%
\mathfrak{h}^{\ast }$, denoted $W$ and called \emph{the Weyl group of} $%
\mathfrak{h}$ \emph{in} $\mathfrak{g}$.

Let $Z(\mathfrak{g})$ denote the center of the enveloping algebra $U(%
\mathfrak{g})$ of $\mathfrak{g}$. \emph{A} $\mathfrak{g}$-\emph{%
infinitesimal character} $\Theta $ is a homomorphism of algebras 
\begin{equation*}
\Theta :Z(\mathfrak{g})\rightarrow \mathbb{C}\text{.}
\end{equation*}%
Since $Z(\mathfrak{g})$ acts on an irreducible Harish-Chandra module (and
also any corresponding smooth globalization) by a scalar, the infinitesimal
character is an important invariant associated to an irreducible, admissible
representation. We now recall Harish-Chandra's parametrization of
infinitesimal characters. We choose a Borel subalgebra $\mathfrak{b}$
containing $\mathfrak{h}$. Thus 
\begin{equation*}
\mathfrak{b}=\mathfrak{h}\oplus \mathfrak{n}\text{ \ where \ }\mathfrak{n}=%
\left[ \mathfrak{b},\mathfrak{b}\right] \text{ is the nilradical of }%
\mathfrak{b}\text{.}
\end{equation*}%
Then one knows that $Z(\mathfrak{g})\subseteq U(\mathfrak{h})\oplus U(%
\mathfrak{g})\mathfrak{n}$ and that the corresponding projection of $Z(%
\mathfrak{g})$ in $U(\mathfrak{h})$ defines an injective morphism of
algebras called \emph{the unnormalized Harish-Chandra map}. We can use this
morphism to identify infinitesimal characters with Weyl group orbits in $%
\mathfrak{h}^{\ast }$ in the following way. Let $\rho $ denote one-half the
sum of the positive roots and suppose $\lambda \in \mathfrak{h}^{\ast }$. \
Then, via the unnormalized Harish-Chandra map, the composition 
\begin{equation*}
\Theta :Z(\mathfrak{g})\rightarrow U(\mathfrak{h})\overset{\lambda +\rho }{%
\rightarrow }\mathbb{C}
\end{equation*}%
defines an infinitesimal character $\Theta $. One knows that for $w\in W$,
the element $w\lambda \in \mathfrak{h}^{\ast }$ defines the same
infinitesimal character $\Theta $. Abusing notation somewhat, we write $%
\Theta =W\cdot \lambda $. The infinitesimal character is called \emph{regular%
} when the only element of $W$ fixing an element in the orbit $W\cdot
\lambda $, is the identity. This is equivalent to the condition that 
\begin{equation*}
\overset{\vee }{\alpha }(\lambda )\neq 0\text{ \ for each }\alpha \in \Sigma 
\text{.}
\end{equation*}

For a $\mathfrak{g}$-module $M$ with regular infinitesimal character one has
the following result. The notes by D. Milicic \cite{M} contain a proof.

\begin{theorem}
Let $M$ be a $\mathfrak{g}$-module with regular infinitesimal character $%
\Theta $. Suppose $\mathfrak{b}$ is a Borel subalgebra of $\mathfrak{g}$
with Levi decomposition 
\begin{equation*}
\mathfrak{b}=\mathfrak{h}\oplus \mathfrak{n}.
\end{equation*}%
Let $\lambda \in \mathfrak{h}^{\ast }$ such that $\Theta =W\cdot \lambda $
and let $\rho $ be one half the sum of the positive roots. Then the Cartan
subalgebra $\mathfrak{h}$ acts semisimply on the $\mathfrak{n}$-homology
groups $H_{\text{p}}(\mathfrak{n},M)$ with eigenvalues of the form $w\lambda
+\rho $ for $w\in W$. In particular 
\begin{equation*}
H_{\text{p}}(\mathfrak{n},M)=\dbigoplus\limits_{w\in W}H_{\text{p}}(%
\mathfrak{n},M)_{w\lambda +\rho }
\end{equation*}%
where 
\begin{equation*}
H_{\text{p}}(\mathfrak{n},M)_{w\lambda +\rho }=\left\{ v\in H_{\text{p}}(%
\mathfrak{n},M):\xi \cdot v=\left( w\lambda +\rho \right) (\xi )v\text{ for
each }\xi \in \mathfrak{h}\right\} .
\end{equation*}
\end{theorem}

A generalization of this result works for any parabolic subalgebra $%
\mathfrak{p}$ with Levi decomposition 
\begin{equation*}
\mathfrak{p}=\mathfrak{l}\oplus \mathfrak{n.}
\end{equation*}%
In particular, if $M$ is a $\mathfrak{g}$-module $M$ with regular
infinitesimal character $\Theta $ and if $Z(\mathfrak{l})$ denotes the
center of the enveloping algebra of $\mathfrak{l}$ then $H_{\text{p}}(%
\mathfrak{n},M)$ is a semisimple $Z(\mathfrak{l})$-module and decomposes
into a direct sum of $Z(\mathfrak{l})$-eigenspaces, where the associated $%
\mathfrak{l}$-infinitesimal characters that appear are related to $\Theta $
by an appropriately defined Harish-Chandra map.

\section{Kostant's theorem}

When $G_{0}$ is a connected, compact Lie group, there is a result, called
Kostant's theorem, that calculates the $\mathfrak{n}$-homolgy groups of an
irreducible representation. In this section we review that result, with
special emphasis on the case of a Borel subalgebra.

Assume $G_{0}$ is a compact real form of $G$. Fix a Borel subalgebra $%
\mathfrak{b}$. Then the normalizer of $\mathfrak{b}$ in $G_{0}$ is a maximal
torus $H_{0}$ and a real form for a Cartan subgroup $H$ of $G$. We let $%
\mathfrak{h}$ be the Lie algebra of $H$ and $\mathfrak{n}$ the nilradical of 
$\mathfrak{b}$. $\Sigma $ is the set of roots. The roots of $\mathfrak{h}$
in $\mathfrak{b}$ determine a set of positive roots $\Sigma ^{+}\subseteq
\Sigma $. Let $\rho $ be one-half the sum of the positive roots. Suppose 
\begin{equation*}
\chi :H_{0}\rightarrow \mathbb{C}
\end{equation*}%
is a continuous character and let $\mu \in \mathfrak{h}^{\ast }$ denote the
complexification of the derivative of $\chi $. To be consistent with the
notation in Section 7 we use \emph{the shifted parameter} 
\begin{equation*}
\lambda =\mu -\rho \text{.}
\end{equation*}%
One knows that 
\begin{equation*}
\overset{\vee }{\alpha }(\lambda )\text{ \ is an integer for each }\alpha
\in \Sigma \text{.}
\end{equation*}%
The character $\chi $ is called \emph{antidominant and regular} if 
\begin{equation*}
\overset{\vee }{\alpha }(\lambda )\notin \left\{ 0,1,2,3,\ldots \right\} 
\text{ \ for each }\alpha \in \Sigma ^{+}.
\end{equation*}%
The Cartan-Weyl parametrization of irreducible representations is as follows.

\begin{theorem}
Maintain the established notations. \newline
\textbf{(a)} Suppose $M$ is an irreducible $G_{0}$-module. Then the space of
coinvariants $H_{\text{0}}(\mathfrak{n},M)$ is an irreducible $H_{0}$-module
and the associated character $\chi $ is antidominant and regular. This
character is called the lowest weight. \newline
\textbf{(b)} If two irreducible representations have the same lowest weight
then they are isomorphic. \newline
\textbf{(c)} To each antidominant and regular character there is an
irreducible $G_{0}$-module with the given character as its lowest weight.
\end{theorem}

We need to define the length function on the Weyl group. One knows that the
Weyl group permutes the roots of $\mathfrak{h}$ in $\mathfrak{g}$. We can
define \emph{the length of} $w\in W$ to be the number of roots in 
\begin{equation*}
-\Sigma ^{+}\cap w\Sigma ^{+}\text{.}
\end{equation*}%
Kostant's theorem is the following:

\begin{theorem}
Suppose $M$ is the irreducible representation for $G_{0}$ with lowest weight 
$\chi $ and let $\lambda \in \mathfrak{h}^{\ast }$ be the shifted parameter.
Then $H_{\text{p}}(\mathfrak{n},M)$ is a sum of irreducible $H_{0}$-modules
each having multiplicity one. The characters of $H_{0}$ that show up as
eigenvalues in $H_{\text{p}}(\mathfrak{n},M)$ are exactly those whose
derivative have the form $w\lambda +\rho $ where the length of $w$ is p$.$
\end{theorem}

In the more general case of a parabolic subalgebra $\mathfrak{p}$, let $%
L_{0} $ be the normalizer of $\mathfrak{p}$ in $G_{0}$ and let $\mathfrak{l}$
be the complexified Lie algebra of $L_{0}$. One knows that $L_{0}$ is
connected and that $\mathfrak{l}$ is a Levi factor of $\mathfrak{p}$.
Indeed, if $L$ is the connected subgroup of $G$ with Lie algebra $l$ then $%
L_{0}$ is the compact real form of $L$. Suppose $M$ is an irreducible
representation for $G_{0}$ and let $\mathfrak{n}$ be the nilradical of $%
\mathfrak{p}$. Then Kostant's Theorem describes the structure of the pth
homology group $H_{\text{p}}(\mathfrak{n},M)$ as an $L_{0}$-module. In
particular, the theorem states that an irreducible representation $V$ of $%
L_{0}$ has, at most, multiplicity one in $H_{\text{p}}(\mathfrak{n},M)$ and
gives a precise condition when $V$ appears, in terms of the degree p, the
lowest weight of $M $ and the lowest weight of $V$. We refer the reader to 
\cite[Chapter IV, Section 9]{KV} for more details.

\section{Flag manifolds and comparison theorems}

As we mentioned before, when $G_{0}$ is noncompact, there are several
conjugacy classes of Borel subalgebras and the structure of the $\mathfrak{n}
$-homology groups of a representation can depend on \ the choice of $G_{0}$%
-conjugacy class. On the other hand, when $M$ is a Harish-Chandra module,
then the locally finite $K_{0}$-action on $M$ extends naturally to a locally
holomorphic $K$-action, and it turns out that the $\mathfrak{n}$-homology
groups of $M$ depend on the $K$-conjugacy classes of Borel subalgebras. In
order to compare the $\mathfrak{n}$-homology groups of $M$ with the $%
\mathfrak{n}$-homology groups of a smooth globalization, we therefore need
to know something about the relationship between $G_{0}$-conjugacy classes
and $K$-conjugacy classes. There is an elegant geometric result, referred to
as \emph{Matsuki duality}, that gives us the needed information. We now
review that result.

One knows that the group $G$ acts transitively on the set of Borel
subalgebras of $\mathfrak{g}$. The corresponding $G$-homogeneous complex
manifold $X$ is called \emph{the full flag space}. In general, if $\mathfrak{%
p}$ is a parabolic subalgebra of $\mathfrak{g}$ then the normalizer of $%
\mathfrak{p}$ in $G$ is the connected subgroup $P$ with Lie algebra $%
\mathfrak{p}$ and the corresponding quotient 
\begin{equation*}
Y=G/P
\end{equation*}%
is called \emph{a flag manifold}. The points in $Y$ are naturally identified
with the $G$-conjugates to $\mathfrak{p}$.

Let $\theta :\mathfrak{g}\rightarrow \mathfrak{g}$ be the complexification
of a Cartan involution of $\mathfrak{g}_{0}$ corresponding to the maximal
compact subgroup $K_{0}$. A Cartan subalgebra $\mathfrak{h}$ of $\mathfrak{g}
$ is called \emph{stable} if $\mathfrak{g}_{0}\cap \mathfrak{h}$ is a real
form and if $\theta (\mathfrak{h})=\mathfrak{h}$. A Borel subalgebra is
called \emph{very special} if it contains a stable Cartan subalgebra. A
stable Cartan subalgebra of a Borel subalgebra is unique (when it exists). A
point in the full flag space is called \emph{very special} if the
corresponding Borel subalgebra is.

Matsuki has established the following \cite{M1}.

\begin{theorem}
Let $X$ be the full flag space. Then \newline
\textbf{(a)} The subset of very special points in a $G_{0}$-orbit is a
nonempty $K_{0}$-orbit. \newline
\textbf{(b)} The subset of very special points in a $K$-orbit is a nonempty $%
K_{0}$-orbit.
\end{theorem}

It follows that the very special points give a one-to-one correspondence
between the $G_{0}$-orbits and the $K$-orbits on $X$, defined by the
following duality. A $G_{0}$-orbit $S$ is said to be \emph{dual} to a $K$%
-orbit $Q$ when $S\cap Q$ contains a special point. In this duality, open $%
G_{0}$-orbits correspond to closed $K$-orbits and the (unique) closed $G_{0}$%
-orbit corresponds to the (unique) open $K$-orbit. We note that Matsuki has
established a similar result for any flag manifold \cite{M2}.

\begin{example}
Suppose $G=SL(2,\mathbb{C)}$, the group of $2\times 2$ complex matrices with
determinant 1 and let $G_{0}=SL(2,\mathbb{R})$. Then the full flag space $X$
is isomorphic to the Riemann sphere. $G_{0}$ has three orbits on $X$. The
closed $G_{0}$-orbit can be identified with an equatorial circle and the
other two orbits are the corresponding open hemispheres. It turns out every
point in the closed orbit is very special, independent of the choice of $%
K_{0}$ (this is true in general for the closed orbit). Put $K_{0}=SO(2,%
\mathbb{R})$. Thus $K=SO(2,\mathbb{C})$. Then the three $K$-orbits on $X$
are a punctured plane, containing the closed $G_{0}$-orbit, and two fixed
points, which can be identified with the respective poles in each of the
open hemispheres. These two poles are the other very special points.
\end{example}

When $M$ is a Harish-Chandra module and $\mathfrak{n}$ is the nilradical of
a Borel subalgebra then one knows that the homology groups $H_{\text{p}}(%
\mathfrak{n},M)$ are finite-dimensional, so it may seem reasonable to ask
when $H_{\text{p}}(\mathfrak{n},M)$ coincides with the $\mathfrak{n}$%
-homolgy groups of a smooth globalization. It turns out this not only
depends on the choice of Borel subalgebra, but also in the the choice of
smooth globalization. When $\mathfrak{n}$ is the nilradical of a very
special Borel subalgebra, $M$ is a Harish-Chandra module, and $M_{\text{min}%
} $ is the minimal globalization, then H. Hecht and J. Taylor have shown 
\cite{HT2} that the natural map 
\begin{equation*}
M\rightarrow M_{\text{min}}\text{ \ induces isomorphisms }H_{\text{p}}(%
\mathfrak{n},M)\rightarrow H_{\text{p}}(\mathfrak{n},M_{\text{min}}).
\end{equation*}%
On the other hand, for the maximal globalization, there are counterexamples
to this result.

The result of Hecht and Taylor has been generalized in the following form. A
Levi factor $\mathfrak{l}$ of a parabolic subalgebra $\mathfrak{p}$ is
called \emph{stable} if $\mathfrak{l}\cap \mathfrak{g}_{0}$ is a real form $%
\mathfrak{l}_{0}$ of $\mathfrak{l}$ and if $\theta (\mathfrak{l})=\mathfrak{l%
}$. The parabolic subalgebra $\mathfrak{p}$ is called \emph{very special} if
it contains a stable Levi factor. Such a Levi factor is unique. Unlike the
case of the full flag space, there may be parabolic subalgebras which are
not $G_{0}$-conjugate to a very special parabolic subalgebra, so we are not
considering all orbits on every flag manifold. However, suppose $\mathfrak{p}
$ is very special and $\mathfrak{l}$ is the stable Levi factor. Define $%
L_{0} $ to be the subgroup of $G_{0}$ that normalizes $\mathfrak{p}$ and
normalizes $\mathfrak{l}$. Then $L_{0}$ is a linear reductive Lie group with
complexified Lie algebra $\mathfrak{l}$ and maximal compact subgroup $%
L_{0}\cap K_{0}$, called \emph{the associated real Levi subgroup}. We have
have the following result \cite[Proposition 2.24]{HS}.

\begin{proposition}
Suppose $\mathfrak{p}$ is a very special parabolic subalgebra with $L_{0}$
and $\mathfrak{l}$ defined as above. Let $\mathfrak{n}$ be the nilradical of 
$\mathfrak{p}$ and suppose $M$ is a Harish-Chandra module for $(\mathfrak{g}%
,K_{0})$. Then the $\mathfrak{n}$-homology groups are Harish-Chandra modules
for $(\mathfrak{l}$,$K_{0}\cap L_{0})$.
\end{proposition}

For the minimal globalization, we have the following \cite{Br2}.

\begin{theorem}
Maintain the hypothesis of the previous proposition. Then the standard
complex induces a Hausdorff topology on $H_{\text{p}}(\mathfrak{n},M_{\text{%
min}})$ and the natural map $M\rightarrow M_{\text{min}}$ induces
isomorphisms 
\begin{equation*}
H_{\text{p}}(\mathfrak{n},M)_{\text{min}}\cong H_{\text{p}}(\mathfrak{n},M_{%
\text{min}})\text{.}
\end{equation*}
\end{theorem}

One might conjecture that above theorem works for the smooth globalization,
and W. Casselman has informed the author that he has proven something along
these lines, although details are unclear. Two partial comparison theorems
about smooth globalizations have been published by other mathematicians. H.
Hecht and J. Taylor have shown the result for minimal parabolic subgroups of 
$G_{0}$ \cite{H}, while U. Bunke and M. Olbrich have shown the result for
any real parabolic subgroup \cite{B}.

D. Vogan has conjectured that all four canonical globalizations commute with
the $\mathfrak{n}$-homology groups of a very special parabolic subalgebra
when the corresponding $G_{0}$-orbit on the flag manifold is open \cite{V1}.
We remark that it has recently been shown that Vogan's conjecture is true
for one globalization if and only it's true for the dual \cite{Br3}. Thus
the conjecture is proven for both the minimal and maximal globalization.

\section{The $\mathfrak{n}$-homology of standard modules}

In the noncompact case, the problem of calculating $\mathfrak{n}$-homology
groups can be quite complicated and there seems to be little hope of just
writing down a formula that generalizes Kostant's theorem for all
irreducible representations. However, there are certain representations,
called \emph{standard modules}, whose $\mathfrak{n}$-homology groups are a
bit more predictable. These standard modules are generically irreducible,
coincide with irreducibles when $G_{0}$ is compact, and can be used to
classify the irreducible representations. In this section we define the
standard representations and consider their $\mathfrak{n}$-homology groups,
focusing on the case of the full flag space.

In particular, we use the construction of minimal globalizations given in 
\cite{HT1}. Let $X$ be the full flag space and, since we need to keep track
of points, introduce the following notation. For $x\in X$ we let $\mathfrak{b%
}_{x}$ be the corresponding Borel subalgebra and let $\mathfrak{n}_{x}$
denote the nilradical of $\mathfrak{b}_{x}$. When we are interested in
calculating the $\mathfrak{n}_{x}$-homology of Harish-Chandra modules, we
can assume $\mathfrak{b}_{x}$ is a very special Borel subalgebra. In that
case, $\mathfrak{h}_{x}$ denotes the stable Cartan subalgebra of $\mathfrak{b%
}_{x}$ and $H_{0}\subseteq G_{0}$ is the corresponding real Cartan subgroup
(thus $H_{0}$ is the associated Levi subgroup). By our linear assumptions on 
$G_{0}$, it follows that $H_{0}$ is abelian, so that an irreducible,
admissible representation of $H_{0}$ is a continuous character 
\begin{equation*}
\chi :H_{0}\rightarrow \mathbb{C}\text{.}
\end{equation*}%
Let $S\subseteq X$ be the $G_{0}$-orbit of $x$. In a natural way, $\chi $
extends to a character of the normalizer of $\mathfrak{b}_{x}$ in $G_{0}$
(we note that $H_{0}$ and the normalizer of $\mathfrak{b}_{x}$ coincide
exactly when $S$ is open). Thus $\chi $ determines a $G_{0}$-homogeneous
analytic line bundle over $S$. One can then define the concept of \emph{a
polarized section \cite[Section 8]{HT1}}. When $S$ is open, the polarized
sections are holomorphic sections, and in general the polarized sections are
locally isomorphic with the restricted holomorphic functions. Let $\mathcal{A%
}(x,\chi )$ denote the sheaf of polarized sections on $S$ and, for p$%
=0,1,2,3,\ldots $ let 
\begin{equation*}
H_{\text{c}}^{\text{p}}(S,\mathcal{A}(x,\chi ))
\end{equation*}%
denote the corresponding compactly supported sheaf cohomology group. Suppose 
$\mu \in $ $\mathfrak{h}_{x}^{\ast }$ is the complexified differential of $%
\chi $, let $\rho $ be one-half the sum of the positive roots for $\mathfrak{%
h}_{x}$ in $\mathfrak{b}_{x}$ and let $\lambda \in \mathfrak{h}_{x}^{\ast }$
denote the corresponding shifted parameter. Thus 
\begin{equation*}
\mu =\lambda +\rho .
\end{equation*}%
We have the following theorem \cite{HT1}.

\begin{theorem}
Maintain the previously defined notations. Let q be the codimension of the $%
K $-orbit of $x$ in $X$. \newline
\textbf{(a)} $H_{\text{c}}^{\text{p}}(S,\mathcal{A}(x,\chi ))$ carries a
natural topology and a continuous $G_{0}$-action, so that the resulting
representation is a minimal globalization. \newline
\textbf{(b)} $H_{\text{c}}^{\text{p}}(S,\mathcal{A}(x,\chi ))$ has
infinitesimal character $\Theta =W\cdot \lambda .$. \newline
\textbf{(c)} When $\lambda $ is antidominant then $H_{\text{c}}^{\text{p}}(S,%
\mathcal{A}(x,\chi ))=0$ when p$\neq $q. \newline
\textbf{(d)} When $\lambda $ is antidominant and $\Theta $ is regular then $%
H_{\text{c}}^{\text{q}}(S,\mathcal{A}(x,\chi ))$ contains a unique
irreducible submodule. In particular, $H_{\text{c}}^{\text{q}}(S,\mathcal{A}%
(x,\chi ))\neq 0$.
\end{theorem}

When $\lambda $ is antidominant and $\Theta $ is regular, we call $H_{\text{c%
}}^{\text{q}}(S,\mathcal{A}(x,\chi ))$ \emph{a regular standard module.}
These modules can be used to parametrize irreducible representations with
regular infinitesimal character. For the remainder of this article we will
make some remarks about how to calculate the $\mathfrak{n}$-homology of
regular standard modules. But we first note that, in the case of a singular
infinitesimal character, the definition of standard module is more subtle,
and the calculation of $\mathfrak{n}$-homology is more elusive.

To state results we will need to differentiate points where we calculate $%
\mathfrak{n}$-homology and the corresponding parameters for eigenvalues of a
Cartan subalgebra (see Theorem 4.1). In particular, we fix a very special
point $x\in X$ as a base point. For $\lambda \in \mathfrak{h}_{x}^{\ast }$
we put $\lambda (x)=\lambda $. When $\mathfrak{b}_{y}$ is a very special
Borel subalgebra and $\mathfrak{h}_{y}\subseteq \mathfrak{b}_{y}$ is the
stable Cartan subalgebra, then there exists $g\in G$ such that 
\begin{equation*}
g\mathfrak{b}_{x}g^{-1}=\mathfrak{b}_{y}\text{ \ and \ }g\mathfrak{h}%
_{x}g^{-1}=\mathfrak{h}_{y}\text{.}
\end{equation*}%
Thus 
\begin{equation*}
\left( \mathfrak{h}_{x}^{\ast }\right) ^{g}=\mathfrak{h}_{y}
\end{equation*}%
This isomorphism is independent of the choice of $g\in G$. For $\lambda \in 
\mathfrak{h}_{x}^{\ast }$ put 
\begin{equation*}
\lambda (y)=\lambda ^{g}\in \mathfrak{h}_{y}^{\ast }.
\end{equation*}%
We note that $\alpha \in \mathfrak{h}_{x}^{\ast }$ is a root of $\mathfrak{h}%
_{x}$ in $\mathfrak{g}$ $\Leftrightarrow $ $\alpha (y)$ is a root of $%
\mathfrak{h}_{y}$ and that $\alpha $ is positive at $x$ $\Leftrightarrow $ $%
\alpha (y)$ is positive at $y$. In particular, $\rho (y)$ is one-half the
sum of the positive roots for $\mathfrak{h}_{y}$ in $\mathfrak{b}_{y}$.

The circle of ideas utilized in \cite{HT1} depend on an identification of
the derived functor of $\mathfrak{n}$-homolgy, in a certain weight (Theorem
4.1), with the geometric fiber applied to a certain, corresponding
localization functor. These ideas originate in the an elegant generalization
of Casselman's submodule theorem, given by A. Beilinson and J. Bernstein in 
\cite{BB}. This identification, together with some functorial rigmarole,
immediately leads to the following result.

\begin{proposition}
Suppose $V=H_{\text{c}}^{\text{q}}(S,\mathcal{A}(x,\chi ))$ is a regular
standard module. Maintain the previously introduced notations. Let $\mathbb{C%
}_{\chi }$ denote the 1-dimensional representation of $H_{0}$ corresponding
to $\chi $ and let $\mathfrak{n}_{x}$ be the nilradical of $\mathfrak{b}_{x}$%
. Then we have the following. \newline
\textbf{(a)} 
\begin{equation*}
H_{\text{p}}(\mathfrak{n}_{x},V)_{\lambda +\rho }=0\text{ \ for p}\neq \text{%
q and }H_{\text{q}}(\mathfrak{n}_{x},V)_{\lambda +\rho }=\mathbb{C}_{\chi }%
\text{. \ }
\end{equation*}%
\newline
\textbf{(b)} If $\mathfrak{b}_{y}$ is a very special Borel subalgebra with
nilradical$\mathfrak{\ n}_{y}$ and $y\notin S$ then 
\begin{equation*}
H_{\text{p}}(\mathfrak{n}_{y},V)_{\lambda (y)+\rho (y)}=0\text{ \ for each p.%
}
\end{equation*}
\end{proposition}

According to Theorem 4.1, the problem of calculating of the $\mathfrak{n}%
_{y} $-homology groups of $V$, at a special point $y\in X$, reduces to the
problem of calculating the values in the weights $\left( w\lambda +\rho
\right) (y)$ for $w\in W$. In geometric terms, this means calculating the
geometric fibers of certain localizations of $V$ or, equivalently,
calculating the result of the so called\emph{\ intertwining functor}. We
briefly consider this problem.

In general, a positive root is called \emph{simple} if it cannot be
decomposed into a nontrivial sum of positive roots. Let $\Sigma _{y}^{+}$ be
the positive roots associated to $\mathfrak{h}_{y}$ in $\mathfrak{b}_{y}$.
For a simple root $\alpha \in \Sigma _{x}^{+}$, the problem of calculating
the values of the $\mathfrak{n}_{y}$-homology groups in the weight $\left(
s_{\alpha }\lambda +\rho \right) (y)$ can be geometrically reduced to
specific calculations for certain real subgroups of $SL(2,\mathbb{C})$. To a
large extent, this idea is already exploited and explained in \cite{HT1} and
some of the necessary calculations are dealt with there.

We finish with an example where, using these ideas, a general formula, like
Kostant's, can be actually written down. Assume $G_{0}$ is a connected,
complex reductive Lie group. Fix a very special Borel subalgebra $\mathfrak{b%
}_{x}$ with stable Cartan subalgebra $\mathfrak{h}_{x}$ and let $\Sigma
_{x}^{+}$ denote the corresponding positive roots. For each $w\in W$, the
set 
\begin{equation*}
w(\Sigma _{x}^{+})=\Sigma _{w\cdot x}^{+}
\end{equation*}%
defines a new set of positive roots and thus a corresponding Borel
subalgebra $\mathfrak{b}_{w\cdot x}$ of $\mathfrak{g}$ containing the stable
Cartan subalgebra $\mathfrak{h}_{x}$. Thus the point $w\cdot x\in X$ is very
special. Because $G_{0}$ is a complex reductive group, one knows that each
Borel subalgebra of $\mathfrak{g}$ is $G_{0}$-conjugate to a Borel
subalgebra of the form $\mathfrak{b}_{w\cdot x}$. Suppose $H$ is the Cartan
subgroup of $G$ with Lie algebra $\mathfrak{h}_{x}$. Then each $\alpha \in
\Sigma $ defines a holomorphic character of $H$ and by restriction, a
corresponding character of $H_{0}$. We write $\chi _{\alpha }$ for this
character of $H_{0}$. If we let $\lambda $ be the shifted parameter and put $%
y=w\cdot x$ then $\left( w\lambda \right) (y)=\lambda $.

Using the above ideas, one can deduce the following.

\begin{theorem}
Let $G_{0}$ be a connected, complex reductive group. Suppose $V=H_{\text{c}%
}^{\text{q}}(S,\mathcal{A}(x,\chi ))$ is the previously defined regular
standard module and assume the $G_{0}$-orbit of $x$ is open in $X$. We
define a chain $\alpha _{1},\ldots ,\alpha _{k}$ of simple roots to be a
finite sequence of roots of $\mathfrak{h}_{x}$ such that for each $j$, $%
\alpha _{j+1}$ is simple for the set of positive roots defined by 
\begin{equation*}
s_{\alpha _{1}}s_{\alpha _{2}}\cdots s_{\alpha _{j}}(\Sigma _{x}^{+})\text{. 
}
\end{equation*}%
Suppose $\alpha _{1},\ldots ,\alpha _{k}$ is a chain of simple roots and let 
$w\in W$ be the ordered product of reflections given by this chain. Let $%
\chi _{w}$ be the character of $H_{0}$ defined by 
\begin{equation*}
\chi _{w}=\chi _{\alpha _{1}}^{-1}\chi _{\alpha _{2}}^{-1}\cdots \chi
_{\alpha _{k}}^{-1}
\end{equation*}%
and let $\mathbb{C}_{\chi }\otimes \mathbb{C}_{\chi _{w}}$ be the
1-dimensional representation of $H_{0}$ corresponding to the character $\chi
\cdot \chi _{w}$. Let q$^{w}$ denote the codimension of the $K$-orbit of $%
y=w\cdot x$ in $X$. Then 
\begin{equation*}
H_{\text{p}}(\mathfrak{n}_{y},V)=0\text{ \ for p}\neq \text{q}^{w}\text{ and 
}H_{\text{q}^{w}}(\mathfrak{n}_{y},V)=\mathbb{C}_{\chi }\otimes \mathbb{C}%
_{\chi _{w}}\text{. \ }
\end{equation*}%
\newline
\end{theorem}

We note that the hypothesis of the theorem implies that the representation $%
V $ is irreducible, and also remark that any attempt to write down a similar
result for other orbits (even in the case of a connected, complex reductive
group), when the standard module is reducible, is considerably more
complicated.

\end{document}